\documentclass{amsart}
\usepackage{amsthm,amsfonts,amsmath,amscd,amssymb,latexsym,epsfig}

\newcommand{\su}{{\mathfrak s  \mathfrak u}}

\newcommand{\g}{{\mathfrak g}}         

\newcommand{\cx}{{\mathbb C}}

\newcommand{\ad}{\operatorname{ad}}
\newcommand{\Ad}{\operatorname{Ad}}

\newcommand{\im}{\operatorname{Im}}

\newcommand{\Lie}{\operatorname{Lie}}

\newcommand{\Ker}{\operatorname{Ker}}

\numberwithin{equation}{section}

\newtheorem{theorem}{Theorem}[section]

\newtheorem{lemma}[theorem]{Lemma}

\newtheorem{corollary}[theorem]{Corollary}

\newtheorem{proposition}[theorem]{Proposition}

\theoremstyle{remark}

\newtheorem{remarknon}{Remark}[section]

\newtheorem{definition}{Definition}[section]

\newtheorem{ack}{Acknowledgment}

\newcommand{\oC}{{\mathbb{C}}}

\newcommand{\oH}{{\mathbb{H}}}

\newcommand{\oP}{{\mathbb{P}}}

\newcommand{\oR}{{\mathbb{R}}}

\newcommand{\sA}{{\mathcal{A}}}   

\newcommand{\sG}{{\mathcal{G}}}   

\newcommand{\sM}{{\mathcal{M}}}   

\newcommand{\fG}{{\mathfrak{g}}}

\newcommand{\fP}{{\mathfrak{p}}}

\newcommand{\fS}{{\mathfrak{s}}}

\begin{document}

\title{Curvature of hyperk\"ahler quotients}
\author{Roger Bielawski}
\address{School of Mathematics\\ University of Leeds\\Leeds LS2 9JT, UK}
\address{Mathematisches Institut, Universit\"at G\"ottingen, G\"ottingen 37073, Germany} 



\begin{abstract} We prove estimates for the sectional curvature of hyperk\"ahler quotients and give applications to moduli spaces of solutions to Nahm's equations and Hitchin's equations.
\end{abstract}

\maketitle

\thispagestyle{empty}
This note was motivated by the following observation: the sectional curvature of the moduli space of charge $k$ $SU(2)$-monopoles is bounded (by an
explicit constant depending on normalisations). Unlike most statements about monopole metrics, this one has a remarkably easy proof
which led us to investigate estimates on sectional curvature of general (finite or infinite-dimensional) K\"ahler and hyperk\"ahler
quotients.

We recall that there is an explicit formula, due to J. Jost and X.-W. Peng \cite{JP}, for the sectional curvature of a large class of
quotients, which include hyperk\"ahler quotients. The quotients in \cite{JP} are formed by taking a Riemannian Banach manifold $(M,g)$ with
a smooth and isometric action of Banach Lie group $G$ which is free on an invariant the level set $\phi^{-1}(c)$ of a suitable smooth map
$\phi$. Jost and Peng compute the curvature of $\phi^{-1}(c)/G$ by giving {\em variational} formulae  for the second fundamental form of the embedding $\phi^{-1}(c)\hookrightarrow M$ and for the O'Neill tensor of the submersion
$\phi^{-1}(c)\rightarrow \phi^{-1}(c)/G$.
\par
Our aim is to give only pointwise estimates on the sectional curvature of hyperk\"ahler quotients and so our proofs are much simpler than
in \cite{JP}. We conclude that for $1$- and $2$-dimensional gauge theories, i.e. moduli spaces of solutions to Nahm's equations and to Hitchin's equations, one gets bounds on the curvature for free, i.e. without seeking any apriori bounds on solutions of relevant differential equations. In the $1$-dimensional case, this is a consequence of the Sobolev embedding $W^{1,2}(a,b)\rightarrow L^\infty(a,b)$, while in dimension $2$ this follows from an analogous embedding of $W^{1,2}(Z)$ into the Orlicz space $L_{e^{t^2}-1}(Z)$.
\par
We also give a simple criterion for a hyperk\"ahler quotient of a finite-dimensional vector space to have asymptotically null curvature.

\section{Infinite-dimensional hyperk\"ahler quotients}

\subsection{Riemannian Banach manifolds}
Let $M$ be a smooth Banach manifold modelled on a Banach space $E$ (see \cite{L} for basics on Banach manifolds). We have a well-defined
tangent bundle $TM$ and the cotangent bundle $T^\ast M$ (bundle of continuous linear functionals). Both of these are Banach manifolds.
Since $E^\ast\otimes E^\ast$ is not necessarily complete (with the norm $\|\alpha\|=\sup\{\alpha(x,y); \|x\|_E=\|y\|_E=1\}$), we consider
its completion $E^\ast\widehat{\otimes} E^\ast$ and the corresponding bundle  $T^\ast M \widehat{\otimes} T^\ast M$.
\begin{definition} A {\em weak\,} Riemannian metric on $M$ is a smooth section $g$ of $T^\ast M \widehat{\otimes} T^\ast M$ which induces a (continuous) positive
definite symmetric bilinear form on each tangent space $T_mM$. The metric $g$ is called {\em strong} if the topology induced by $g$ on each fibre is equivalent to
the topology of the model Banach space $E$.\end{definition}

If $(M,g)$ is a weak Riemannian Banach manifold, then the usual proof of existence and uniqueness of the Levi-Civita connection tells us
what $g(\nabla_XY,Z)$ should be for any $Z$ (and so proves the uniqueness), but it does not guarantee existence. However, if  we assume that a
smooth Levi-Civita connection exists, then other Riemannian notions such as parallel transport, geodesics, exponential map, curvature make
sense and have usual properties.
\par
In what follows we shall assume that $(M,g)$ is a weak Riemannian Banach manifold such that the (smooth) Levi-Civita connection exists\footnote{We
resist the temptation to call such manifolds {\em medium} Riemannian Banach manifolds.}.
\par
\begin{definition} A weak hyperk\"ahler Banach manifold is a weak Riemannian manifold $(M,g)$ with Levi-Civita connection $\nabla$ and
three smooth anti-commuting almost complex structures $I_1,I_2,I_3$, which are fibre-wise isometries satisfying $I_1I_2I_3=-1$ and which
commute with $\nabla$. \end{definition} The definition of a weak K\"ahler Banach manifold is analogous.

\subsection{Group actions and quotients} Let $G$ be a Banach Lie group with Lie algebra $\g$. If $G$ acts smoothly on a Banach manifold $M$,
then for any $\rho\in\fG$ we denote by $\check{\rho}$ the corresponding fundamental vector field on $M$. We write $\check{\fG}$ for the ``subbundle"\footnote{In general,$\check{\fG}$ is not locally trivial.}  of the tangent bundle generated by the vector fields $\check{\rho}$.
\par
If $M$ has a weak Riemannian metric $g$, then we define a $\fG^\ast$-valued $1$-form $\Lambda$ by
\begin{equation} \Lambda(v)(\rho)=g(v,\check{\rho}),\label{Ph}
\end{equation}
 i.e.  the pointwise adjoint of the mapping $l_m:\rho\mapsto \check{\rho}_{|_m}$ with respect to $g$.
\begin{definition} A smooth action of a Banach Lie group $G$ on a weak Riemannian manifold $(M,g)$ is called {\em elliptic} at a point $m$ if  $\Lambda(T_mM)=\Lambda(\check{\g}_m)$.\end{definition}
In many infinite-dimensional applications, the map $\rho\mapsto\Lambda(\check{\rho}_{|_m})$ is a second order linear differential operator whose ellipticity guarantees that $\Lambda(T_mM)=\Lambda(\check{\g}_m)$
 \par
 We also observe that the condition of ellipticity at $m$ is equivalent to the addition map $\check{\g}_m\times \check{\g}_m^\perp\rightarrow T_mM$ being an isomorphism (here $\check{\g}_m^\perp$ is the subspace $g$-orthogonal to $\check{\g}_m$). In particular, it implies that $\check{\g}_m$ is a closed subspace.
 \par
Recall that an action is called {\em proper} if the map $G\times M\rightarrow M\times M$, $(g,m)\rightarrow (gm, m)$ is proper.
\begin{proposition} Let there be given a  free, proper, isometric and elliptic action of a Banach Lie group $G$ on a weak Riemannian Banach manifold $(M,g)$ for which a Levi-Civita connection exists. Then the space of orbits $M/G$ is canonically a weak Riemannian Banach manifold with a Levi-Civita connection.\label{quot}\end{proposition}
\begin{proof} Ellipticity and properness  imply that the orbits are closed submanifolds of $M$. One then constructs a slice $S_m$ using the exponential mapping for the Levi-Civita connection at a particular point $m$ of an orbit in the directions of the subbundle $\check{\g}^\perp$. The properness of the action guarantees that $S_m$ can be chosen small enough to be a slice to the action. The properness of the action also implies that $M/G$ is Hausdorff and hence a Banach manifold. The tangent space to $M/G$ at $Gm$ is canonically identified with $\check{\g}^\perp_m$ and this gives us a metric and the Levi-Civita connection.\end{proof}

\subsection{Hyperk\"ahler quotients\label{quotients}}

\begin{definition} A smooth action of $G$ on a weak hyperk\"ahler Banach manifold $(M,g,I_1,I_2,I_3)$ is called tri-Hamiltonian if there exist (moment)
maps $\mu_1,\mu_2,\mu_3:M\rightarrow \g^\ast$ which are smooth, equivariant and satisfy
\begin{equation} <d\mu_i(v),\rho>=g(v,I_i\check{\rho}),\enskip i=1,2,3,\label{moment}\end{equation}
for any tangent vector $v$ and any $\rho\in \g$.\end{definition}

In infinite dimensions, if the metric is only weak, the image of $d\mu_i$ (which is the same as image of $\Lambda$) will only be a dense subspace of $\g^\ast$ and so there is no hope that the moment map $\mu=(\mu_1,\mu_2,\mu_3)$ will be a submersion. We can give a simple criterion for a level set of the moment map to be a manifold. Observe first, that if $M$ is connected and $c_i\in \im \mu_i$, then $\im \mu_i \subset c_i+\im \Lambda$.
\begin{proposition} Let $\mu=(\mu_1,\mu_2,\mu_3)$ be a moment map for a  tri-Hamiltonian action of a Banach Lie group $G$ on a connected weak hyperk\"ahler Banach manifold $M$. Let $c$ be an element of $\g^\ast\otimes \oR^3$ fixed by the coadjoint action of $G$ and such that the action of $G$ is locally free and elliptic at points of $\mu^{-1}(c)$. Suppose that the point-wise image $V=\Lambda(T_mM)$ of $\Lambda$ does not depend on $m$ and that $V$ can be made into a Banach space with respect to a norm $\|\cdot\|_V$, which is   stronger than the one defining the topology of $\g^\ast$ and for which $\mu:M\rightarrow c+V\otimes \oR^3$ remains smooth. Then $\mu^{-1}(c)$ is a submanifold.\label{stronger}\end{proposition}
\begin{proof} Let $M_c=\mu^{-1}(c)$.
Acting by $I_i$ on the splitting $T_mM=\check{\g}_m\oplus \check{\g}_m^\perp$, we have $T_mM=I_i\check{\g}_m\oplus \Ker (d\mu_i)_m$. At a point $m$ of $M_c$, the spaces $I_i\check{g}$ are mutually orthogonal, and, hence, $d\mu_{|_m}$ is an isomorphism between $(\Ker d\mu)^\perp=I_1\check{\g}\oplus I_2\check{\g}\oplus I_3\check{\g}$ and $V \otimes\oR^3$.
Therefore $\mu:M\rightarrow c+V\otimes \oR^3$ is a submersion with respect to a stronger topology on $V$ and, so, $M_c$ is a submanifold.\end{proof}

Finally, we have:
\begin{proposition} Let $\mu=(\mu_1,\mu_2,\mu_3)$ be a moment map for a  tri-Hamiltonian action of a Banach Lie group $G$ on a weak hyperk\"ahler Banach manifold. Let $c$ be an element of $\g^\ast\otimes \oR^3$ fixed by the coadjoint action of $G$ such that $\mu^{-1}(c)$ is a submanifold and the action of $G$ on $\mu^{-1}(c)$  is free, isometric, elliptic and proper.
Then $Q=\mu^{-1}(c)/G$ is a weak hyperk\"ahler Banach manifold.
\label{Q}\end{proposition}
\begin{proof}  Proposition \ref{quot} shows that $Q=M_c/G$ is a Banach manifold. Its tangent space at any $Gm$ is identified with the subspace $H$ of $T_mM$ orthogonal to $\check{\g}\oplus I_1\check{\g}\oplus I_2\check{\g}\oplus I_3\check{\g}$. It is clear that $H$ inherits the Riemannian metric and it is invariant under $I_1,I_2,I_3$. Moreover, the Levi-Civita connection $\overline{\nabla}$
of $Q$ is simply defined by $\overline{\nabla}_XY=\pi_H(\nabla_XY)$, where $\pi_H$ is the orthogonal projection onto $H$ and $X,Y$ are
sections of $H$. The induced complex structures  commute with $\overline{\nabla}$.
\end{proof}

\section{Curvature estimates}

We assume that we are in the situation of Proposition \ref{Q}, i.e. we have a weak hyperk\"ahler manifold
$(M,g,I_1,I_2,I_3)$ with a tri-Hamiltonian action of a Banach Lie group $G$ which is free, proper, isometric and elliptic on the $c$-level
set of the hyperk\"ahler moment map $\mu=(\mu_1,\mu_2,\mu_3)$ ($c\in (\g^\ast\otimes\oR^3)^G$). Thus, at every point $m\in \mu^{-1}(c)$, there is a splitting $T_mM=H\oplus \check{\g}\oplus I_1\check{\g}\oplus I_2\check{\g}\oplus I_3\check{\g} $ as an orthogonal sum of closed subspaces. Moreover, $\mu^{-1}(c)$ is a submanifold and hence the hyperk\"ahler quotient
$Q=\mu^{-1}(c)/G$ exists and is a weak hyperk\"ahler manifold (with Levi-Civita connection) with the tangent bundle of $Q$ identified with $H$. We shall give several estimates on the curvature
of $Q$ in terms of various quantities. We begin  with some definitions. Let $M_c=\mu^{-1}(c)$.
\par
\begin{definition}Let $\check{\g}^\oH=\check{\g}\oplus I_1\check{\g}\oplus I_2\check{\g}\oplus I_3\check{\g} $ and let $m\in M_c$. Denote by $A$ the restriction to
 $\bigl(\check{\g}^\oH\bigr)^\perp\times \bigl(\check{\g}^\oH\bigr)^\perp$ of the O'Neill tensor (cf. \cite{Besse}) of the submersion
 $M\rightarrow M/G$. In other words, for $(X,Y)\in \bigl(\check{\g}^\oH\bigr)_m^\perp\times \bigl(\check{\g}^\oH\bigr)_m^\perp$, $A(X,Y)$
 is the $\check{\g}$-part of $\nabla_{\tilde{X}}\tilde{Y}$ at $m$, where $\nabla$ is the Levi-Civita connection and ${\tilde{X}},\tilde{Y}$ are
 any extensions of $X,Y$ to neighbourhood of $m$.\end{definition}
\begin{definition} We denote by $V:M_c\rightarrow \oR\cup\{+\infty\}$ the pointwise norm of $A$ defined above, i.e.
$$ V(m)=\sup\left\{\frac{A(X,Y)}{|X||Y|}\enskip ; \enskip (X,Y)\in \bigl(\check{\g}^\oH\bigr)^\perp_m\times \bigl(\check{\g}^\oH\bigr)^\perp_m,\enskip X,Y\neq 0\right\}.$$
\end{definition}
Here and in what follows $|\cdot|$ denotes the length of tangent vectors in the metric $g$. Observe that in general the norm of $A$ may be
infinite (since the metric is weak). If, however, the hyperk\"ahler quotient $\mu^{-1}/G$ is finite-dimensional (and hence so is
$\bigl(\check{\g}^\oH\bigr)^\perp$), then $V(m)$ is finite.

We have:
 \begin{proposition} Under the above assumptions the sectional curvature $K_Q$ of the hyperk\"ahler quotient $Q$ of $M$ by $G$ satisfies the pointwise estimate
$$ |K_Q(p)(\pi)-K_M({m})(\tilde{\pi})|\leq 9V(m)^2$$
 where ${m}$ is any point in $M$ projecting to $p$ and $\tilde{\pi}$ is the horizontal lift of a plane $\pi\subset T_pQ$ to $T_{{m}}M$.
\label{main}\end{proposition}
\begin{proof}
Recall that $M_c$ denotes the $c$-level set of the hyperk\"ahler moment map. The space $T_pQ$ is identified with the horizontal subspace of
$T_mM_c$, which in turn is the subspace of $T_mM$ orthogonal to $\check{\g}\oplus I_1\check{\g}\oplus I_2\check{\g}\oplus I_3\check{\g}$.
At points of $M_c$, this last decomposition is orthogonal. Let $X,Y$ be horizontal vector fields on $M_c$ and let us decompose
$\nabla_XY=Z+Z^\perp$ where $Z$ is horizontal and $Z^\perp=\check{\rho}_0+I_1\check{\rho}_1+I_2\check{\rho}_2+I_3\check{\rho}$. Thus
$$\rho_i=A(X,-I_iY)\quad i=0,1,2,3,$$
where $I_0=-1$. Hence
\begin{equation}|\check{\rho}_i|\leq V(m)|Y||X|.\label{fund}\end{equation}

We compute the sectional curvature $K_{M_c}(X,Y)$ of plane in $T_mM_c$ spanned by orthonormal and horizontal vectors $X,Y$. From the Gauss
equation \cite{Besse} (which remains true in infinite-dimensional setting), we know that
$$K_{M_c}(X,Y)=K_M(X,Y)+g\bigl(\alpha(X,X),\alpha(Y,Y)\bigr)-g\bigl(\alpha(X,Y),\alpha(X,Y)\bigr)$$
where $\alpha$ is the second fundamental form of the embedding $M_c\hookrightarrow M$, i.e. $\alpha(X,Y)=\bigl(\nabla_XY\bigr)^\perp$ for
any extension of $Y$ to a vector field near $m$. From the above discussion
$\alpha(X,Y)=I_1\check{\rho}_1+I_2\check{\rho}_2+I_3\check{\rho}_3$ (and similarly for $\alpha(X,X),\alpha(Y,Y)$) and so, using the mutual
orthogonality of $I_i\check{\g}_m$, the estimate \eqref{fund} and the fact $|X|=|Y|=1$, we obtain
\begin{equation}\bigl|K_{M_c}(X,Y)-K_M(X,Y)\bigr|\leq 6 V(m)^2.\label{immersion}\end{equation}
We now compare $K_{M_c}(X,Y)$ to the sectional curvature $K_Q(X,Y)$ in the quotient $Q=M_c/G$. The O'Neill formula \cite{Besse} shows that
\begin{equation}K_Q(X,Y)=K_{M_c}(X,Y)+3\bigl|(\nabla_XY)^v\bigr|^2,\label{submersion}\end{equation}
where the superscript $v$ denotes the vertical
part. In the above notation $(\nabla_XY)^v=\check{\rho}_0$, which together with \eqref{fund} shows that
$$\bigl|K_Q(X,Y)-K_{M_c}(X,Y)\bigr|\leq 3 V(m)^2.$$
\end{proof}

We observe that the proof obviously works as well for K\"ahler quotients and gives a similar estimate. In the K\"ahler case, moreover, the curvature must increase by at least $V(m)$ on some planes:
\begin{proposition} Let $M$ be a weak K\"ahler Banach manifold with a Hamiltonian action  of a Banach Lie group $G$, which is free, proper and strongly isometric on
the $c$-level set of the K\"ahler moment map $\mu$ ($c\in (\g^\ast)^G$).
 The sectional curvature $K_Q$ of the K\"ahler quotient $Q=\mu^{-1}(c)/G$ of $M$ by $G$ satisfies the pointwise estimate
$$ \bigl| K_Q(p)(\pi)-K_M({m})(\tilde{\pi})\bigr|\leq 5V(m)^2$$
 where ${m}$ is any point in $M$ projecting to $p$ and $\tilde{\pi}$ is the horizontal lift of a plane $\pi\subset T_pQ$ to $T_{{m}}M$. Moreover
 $$  \sup_\pi \bigl(K_Q(p)(\pi)-K_M({m})(\tilde{\pi})\bigr)\geq V(m).$$
\label{easy}\end{proposition}
\begin{proof} The proof is essentially the same. The estimate \eqref{immersion} can be replaced by
\begin{equation}\bigl|K_{M_c}(X,Y)-K_M(X,Y)\bigr|\leq 2 V(m)^2.\end{equation}
This, together with \eqref{submersion} and the definition of $V(m)$, proves the estimates.\end{proof}

\medskip

We now wish to give estimates on $V(m)$, which can actually be computed in applications. Let $(M,g,I_1,I_2,I_3)$, $G$ and  $c\in (\g^\ast\otimes \oR^3)^G$ be as above.
Fix an $m\in M_c=\mu^{-1}(c)$. We choose any norm $\|\cdot\|_m$ on $\g$. This norm can be completely different from the norm used to define
the Banach Lie algebra on $\g$. Moreover, this norm may depend on $m\in M_c$.
\begin{definition} Let $B:\bigl(\check{\g}^\oH\bigr)^\perp\times \g\rightarrow \bigl(\check{\g}^\oH\bigr)^\perp$ denote the bilinear map
given by $(X,\rho)\mapsto Z$, where $Z$ is the $\bigl(\check{\g}^\oH\bigr)^\perp$-part of $\nabla_X\check{\rho}$.\end{definition}
We now define $F(m)$ as the norm of $B_m$ with respect to the norm $\|\cdot\|_m$ on $\g$ and $|\cdot|=\sqrt{g(\cdot,\cdot)}$ on $T_mM$:
$$ F(m)=\sup\left\{\frac{B(X,{\rho})}{|X|\cdot\|\rho\|_m}\enskip ; \enskip (X,\rho)\in \bigl(\check{\g}^\oH\bigr)^\perp_m\times \g,\enskip X,\rho\neq 0\right\}.$$

Consider again the $\fG^\ast$-valued $1$-form $\Lambda$ given by $\Lambda(v)(\rho)=g(v,\check{\rho})$. In other words
\begin{equation} \Lambda=-I_1d\mu_1=-I_2d\mu_2=-I_3d\mu_3.\label{Phi}\end{equation}
We denote  by $\|\cdot\|_m^\ast$ the ``norm" induced by $\|\cdot\|_m$ on $\im \Lambda\subset\g^\ast$, i.e. $\|L\|^\ast=\sup
\{L(\rho);\|\rho\|=1\}$. Let $l(m)$ be the norm of $(\Lambda_m)^{-1}:\im \Lambda\rightarrow \check{\g}$:
\begin{equation}l(m)=\sup_{\rho\neq 0}\frac{\bigl|\check{\rho}_m\bigr|}{
\bigl\|\Lambda(\check{\rho}_m)\bigr\|_m^\ast}.\label{l}\end{equation}
Without further assumptions both $l(m)$ and $F(m)$ can be infinite. We have:
\begin{proposition} For any norm $\|\cdot\|_m$ on $\g$ the following inequality holds:
$$ V(m)\leq l(m)F(m).$$
\label{main2}\end{proposition}
\begin{proof} Let $X,Y\in \bigl(\check{\g}^\oH\bigr)^\perp_m$ and let $\tilde{X},\tilde{Y}$ be local sections of  $\bigl(\check{\g}^\oH\bigr)^\perp$
extending $X$ and $Y$. Let $\check{\rho}$ be the fundamental vector field whose value at $m$ is the $\check{\g}$-part of
$\nabla_{\tilde{X}}\tilde{Y}$ (as the action is free at $m$, this is well defined). We need to estimate $|\check{\rho}_m|$. Let $\nu$ be
any element of $\g$.
 Since $\tilde{Y}$ is horizontal, $g(\tilde{Y},\check{\nu})\equiv 0$, and since $\tilde{X}$ is also horizontal
$$ g(\nabla_{\tilde{X}}\tilde{Y},\check{\nu})=-g(\tilde{Y},\nabla_{\tilde{X}}\check{\nu}).$$
At the point $m$ we can rewrite this as
$$\langle \Lambda(\check{\rho}),\nu\rangle=-g(Y,\nabla_X\check{\nu}).$$
Therefore $\|\Lambda(\check{\rho}_m)\|_m^\ast\leq F(m)|X||Y|$ and using the definition of $l(m)$ we get
\begin{equation}|\check{\rho}_m|\leq l(m) F(m)|Y||X|,\quad i=0,1,2,3,\label{fund2}
\end{equation}
which proves the estimate.\end{proof}

Let us discuss this. We start with $F(m)$. From its definition, $F(m)$ is finite if the bilinear operator $B_m$ is continuous for
the norms $\|\cdot\|_m$ and $|\cdot|$. This tells us which norms $\|\cdot\|_m$ are allowed on $\g$. If the hyperk\"ahler quotient is
finite-dimensional, it is easier to decide on $\|\cdot\|_m$:
\begin{lemma} If $\dim Q<+\infty$, then $F(m)$ is bounded providing the linear operator $\rho\rightarrow \bigl(\nabla_X\check{\rho}\bigr)_m$ is
bounded for every $X\in\bigl(\check{\g}^\oH\bigr)_m^\perp$ with respect to the norms $\|\cdot\|_m$ on $\g$ and $|\cdot|$ on
$T_mM$.\end{lemma}
\begin{proof} As $\bigl(\check{\g}^\oH\bigr)^\perp_m$ is finite-dimensional, the bilinear operator $B_m$ is separately continuous in both variables. The
Mazur-Orlicz theorem (see e.g. Corollary $8$ in \cite{Swa}) implies that $B_m$ is continuous.\end{proof}

We observe next that it is easy estimate $F(m)$ for a flat $M$:
\begin{lemma} Let $(M,g)$ be a flat weak hyperk\"ahler manifold, i.e. $(M,g,I_1,I_2,I_3)$ is isomorphic to an open subset of a
quaternionic Banach space $E$ with a continuous bilinear form $g$ on $E$. Suppose that, under this isometry, $\check{\rho}_m=L(\rho,m)
+P(\rho)$, where $L:\g\times E$ is a bilinear operator  and $P$ is independent of $m\in E$. Let $\|\cdot\|$ be a norm on $\g$ such that $L$
is continuous with respect to  $\|\cdot\|$ on $\g$ and $g$ on $E$. Then $F(m)$ is uniformly bounded for $\|\cdot\|_m=\|\cdot\|$ .\hfill $\Box$\label{H}
\end{lemma}
In particular $F$ is bounded for  linear actions of Hilbert Lie groups on Hilbert spaces (with strong metrics).

On the other hand, we have the following estimates on $l(m)$:
\begin{lemma}
$$ l(m)\leq \sup_{\rho\neq 0}\frac{\|\rho\|_m}{\bigl|\check{\rho}_m\bigr|},$$
$$  l(m)\leq \sup_{\rho\neq 0}\left(\frac{\|\rho\|_m} {\bigl\|\Lambda(\check{\rho}_m)\bigr\|_m^\ast}\right)^{1/2}.$$
\label{estimate}\end{lemma}
\begin{proof}
 From the definition of $\Lambda$, we have, for any $\rho\in \g$, $\langle \Lambda(\check{\rho}),\rho\rangle=g(\check{\rho},\check{\rho})$ and hence
$$ |\check{\rho}_m|^2\leq  \|\Lambda(\check{\rho}_m)\|_m^\ast \cdot \|\rho\|_m.$$
The estimates follow thanks to \eqref{l}.
\end{proof}

 Thus $l$ is uniformly bounded (resp. is asymptotically null on $Q$) if, for every $\rho\in\g$ with
$\|\rho\|_m=1$, the length of $\check{\rho}_m$ is bounded away from zero (resp. is asymptotically infinite) on the $c$-level set of the
hyperk\"ahler moment map. In quotients of infinite-dimensional manifolds, it is the second inequality that is useful: the map $\rho\mapsto\Lambda(\check{\rho}_m)$ is often a positive-definite self-adjoint elliptic operator and one easily gets an estimate on $\|\rho\|$ in terms of the norm of $\Lambda(\check{\rho}_m)$.

If we set $\|\rho\|_m=|\check{\rho}_m|$ in the first inequality of the above lemma, we get that $l(m)\leq 1$, and hence we obtain
\begin{corollary} $V(m)$ is bounded by the norm of the second O'Neill tensor at $m$, i.e. the norm of the bilinear operator $C:\bigl(\check{\g}^\oH\bigr)_m^\perp\times \check{\g}_m\rightarrow \bigl(\check{\g}^\oH\bigr)_m^\perp$, $C(X,U)$ is the $\bigl(\check{\g}^\oH\bigr)^\perp$-part of $\bigl(\nabla_X\check{\rho}\bigr)_m$, where $\rho$ is the unique element of  $\g$ such that $\check{\rho}_m=U$. \label{O'Neill}\end{corollary}

\section{Finite-dimensional quotients} We give a simple application to hyperk\"ahler quotients of finite-dimensional vector spaces.
\par
Let $M=\oH^d$ with its Euclidean hyperk\"ahler structure and $G$ be a closed subgroup of  $Sp(d)$ acting linearly on $M$.
If we identify $\fS\fP (d)$ and $\fS\fP (d)^\ast$ with quaternionic matrices $A$ satisfying $A^\dagger=-A$, where $\dagger$ is
transposition followed by quaternionic conjugation, then the hyperk\"ahler  moment map for the action of $Sp(d)$ is
\begin{equation} q\mapsto (qiq^\dagger, qjq^\dagger, qkq^\dagger)\label{dagger}\end{equation}
and we denote by $\mu=(\mu_1,\mu_2,\mu_3)$ the projection of \eqref{dagger} onto three copies of the Lie algebra $\g$ of $G$. $\mu$ is a particular
hyperk\"ahler moment map for $G$.

We have
\begin{theorem} With the above notation, suppose that $G$ acts locally freely on the set $\mu^{-1}(0)-\{0\}$. Let $c=(c_1,c_2,c_3)$ with each $c_i$ a central element of $\g$ and
$G$ acting freely on $\mu^{-1}(c)$. Then the curvature of the hyperk\"ahler quotient $Q=\mu^{-1}(c)/G$ is asymptotically null.\end{theorem}
\begin{proof} We  use Propositions \ref{main} and \ref{main2}. In the finite-dimensional case any norm $\|\cdot\|$ on $\g$ will do.
As observed in Lemma \ref{H}, $F$ is uniformly bounded. We are going to estimate $l$ from the first statement of Lemma \ref{estimate}. Let
$S$ denote the unit sphere in $\oH^d$.  Since $G$ is compact and acts locally freely on $\mu^{-1}(0)-\{0\}$, $G$ acts locally freely on
$V=S\cap \mu^{-1}(D)$, where $D$ is some small closed neighbourhood of $0$ in the center of $\g$. Therefore, for any $\rho\in \g$, the
minimum of $|\check{\rho}|$ over $V$ is a non-zero number, say $\lambda(\rho)$. Consider now the asymptotic behaviour of $|\check{\rho}|$
on $\mu^{-1}(c)$. Let $q\in \mu^{-1}(c)$ and $|q|=R$. Then $q/R\in \mu^{-1}(c/R^2)$ and so, for large enough $R$,
$|\check{\rho}_{\frac{q}{R}}|\geq \lambda(\rho)$. Therefore $|\check{\rho}_{q}|\geq R\lambda(\rho)$ and the result is proven.
\end{proof}

This fact is of course to be expected as $\mu^{-1}(c)/G$ is asymptotically isometric to $\mu^{-1}(0)/G$ and the latter is a cone over a compact orbifold.

The example of  $T^\ast \oP^1\times T^\ast \oP^1$, which is the hyperk\"ahler quotient of $\oH^2\times \oH^2$ by $S^1\times S^1$ (acting
separately on each $\oH^2$), shows that, without the assumption, the quotient does not have to have an asymptotically null curvature.

\section{Moduli spaces of solutions to Nahm's equations}
We wish to estimate the sectional curvature of moduli spaces of solutions to Nahm's equations on an interval with prescribed poles at the end of the interval.
We first describe the hyperk\"ahler quotient construction of this moduli space in order to check that the conditions of Proposition \ref{Q} hold.

\subsection{Construction}
Let $G$ be a compact Lie group with Lie algebra $\g$ and let $\langle\;,\;\rangle$ denote a positive-definite $\Ad G$-invariant inner
product on $\g$ with respect to which
\begin{equation} \bigl| [A,B]\bigr|\leq 2|A||B|, \quad \text{for all $A,B\in \g$}.\label{normalise}\end{equation}
Let $\alpha_i,\beta_i\in \g$, $i=1,2,3$ satisfy $[\alpha_i,\alpha_j]=\epsilon_{ijk}\alpha_k$ and similarly for $\beta_i$, so that these define homomorphisms $\su(2)\rightarrow \g$.
\par
We are going to construct the moduli space of solutions to Nahm's equations on an interval $(a,b)$ with simple poles at $a,b$ and residues $\alpha_i,\beta_i$.
\par
Let $e_1,e_2,e_3$ denote the right multiplication by $i,j,k$ on $\oH$ and consider the linear operators $L_\alpha=\sum_{i=1}^3 (\ad\alpha_i)\otimes e_i$ and  $L_\beta=\sum_{i=1}^3 (\ad\beta_i)\otimes e_i$ on $\g\otimes \oH$.

We define a space $E$ as the space of $\g\otimes\oH$-valued continuously differentiable functions $u$ on $(a,b)$ such that
\begin{equation} L(u)(s)=\frac{du}{ds}-\frac{L_\alpha(u)(s)}{s-a}- \frac{L_\beta(u)(s)}{s-b}\label{L}\end{equation}
is continuous on $[a,b]$. The right multiplication by quaternions  preserves this space. We put a norm on $E$ by
$$ \|u\|=\|u\|_{C^0}+\|L(u)\|_{C^0},$$
where  the $C^0$-norms are the $\sup$-norms defined on $\g\otimes \oH$-valued functions using  the usual norm on the quaternions and the chosen invariant inner product on $\g$. With this norm $E$ is a Banach space (and a closed subspace of $C^1(a,b)$).
\par
We consider the following metric $g$ on $E$: if, after identifying $\oH\simeq \oR^4$, $v=(t_0,t_1,t_2,t_3)$,
$v^\prime=(t_0^\prime,t_1^\prime,t_2^\prime,t_3^\prime)$ are in $E$, then
 \begin{equation}g(v,v^\prime)=\sum_{i=0}^3\int_{a}^b\langle t_i,t^\prime_i \rangle ds \label{metric}.\end{equation}
We now consider an affine space $\sM$ defined as $iS_1(t)+jS_2(t)+kS_3(t)+E$, where $S_i(s)=\frac{\alpha_i}{s-a}+\frac{\beta_i}{s-b}$,
$i=0,1,2,3$. We view $(\sM,g)$ as a flat weak hyperk\"ahler
Banach manifold (modelled on $E$) consisting of $\g$-valued quadruples $(T_0,T_1,T_2,T_3)$ with prescribed boundary behaviour.
\par
We also define $\sG$ as the group of gauge transformations $g:[a,b]\rightarrow G$ whose Lie algebra are maps $\rho:[a,b]\rightarrow \g$  of class $C^2$ on $(a,b)$ satisfying $\rho(a)=\rho(b)=0$ and
\begin{equation} P(\rho)(s)=\frac{d^2\rho}{ds^2}+\sum_{i=1}^3\frac{(\ad^2\alpha_i) \rho(s)}{(s-a)^2}+\sum_{i=1}^3\frac{(\ad^2\beta_i) \rho(s)}{(s-b)^2}\label{P}\end{equation}
continuous on $[a,b]$. We equip this Lie algebra with the Banach norm $\|\rho\|_{C^1}+\|P(\rho)\|_{C^0}$ (all norms are $\sup$-norms on $(a,b)$).
\par
 The group $\sG$ acts smoothly on $\sM$ by
\begin{eqnarray} T_0&\mapsto & \Ad(g)T_0-\dot{g}g^{-1}\nonumber\\ T_i&\mapsto & \Ad(g)T_i\;,\;\;\qquad i=1,2,3.\label{action}\end{eqnarray}
Differentiating, we obtain that for any $\rho\in \Lie \sG$
\begin{equation}\check{\rho}_T=\bigl(-\dot{\rho}+[\rho,T_0],[\rho,T_1],
[\rho,T_2],[\rho,T_3]\bigr).\label{check}\end{equation}
This action is Hamiltonian and the moment map equations are:
\begin{equation}\mu_i(T_0,T_1,T_2,T_3)(\rho)=\int_a^b{\langle\dot{T}_i+[T_0,T_i]-
\frac{1}{2}\sum_{j,k=1,2,3}
\epsilon_{ijk}[T_j,T_k],\rho\rangle}ds.\label{Nahm}\end{equation}
The $0$-level set of $\mu=(\mu_1,\mu_2,\mu_3)$ is given by the Nahm's equations: $\dot{T}_1+[T_0,T_1]-[T_2,T_3]=0$ etc.
We need to check that the conditions of Proposition \ref{Q} are satisfied. The action of $\sG$ is isometric, free and proper everywhere.  The form $\Lambda$ is given by:
\begin{equation}
\Lambda(t_0,t_1,t_2,t_3)(\rho)=\int_a^b{\langle\dot{t}_0+[T_0,t_0]+[T_1,t_1]+[T_2,t_2]+[T_3,t_3],\rho\rangle}ds.\label{perp}
\end{equation}
Thus the image of $\Lambda$ is contained in the subspace $V$ of $(\Lie \sG)^\ast$ given by the pairing $\int_a^b{\langle f,\rho\rangle}ds$ for $f\in C^0([a,b])\otimes \g$. According to Propositions \ref{stronger} and \ref{Q}, we need to show that any element of $V$ can be obtained as $\Lambda(\check{\rho})$. Substituting \eqref{check} into \eqref{perp} expresses $\Lambda(\check{\rho})=\langle f,\cdot\rangle$ as the equation
\begin{equation} -\ddot{\rho}+2[\dot{\rho},T_0]+[\rho,\dot{T}_0]-\sum_{i=0}^3 \bigl[T_i,[T_i,\rho]\bigr]=f.\label{ddot}\end{equation}
It is easy to obtain a solution in $C^2(a,b)$ of this equation (with bounded $C^1(a,b)$-norm), e.g. by finding a solution in $W^{1,1}(a,b)$ using variational methods or by approximation method as in Lemmae 2.18--2.20 in \cite{Don}. It is then automatic that $\rho(t)\in\Lie\sG$.

 Thus, the quotient of the space of solutions to \eqref{Nahm} by $\sG$  is a (finite-dimensional) hyperk\"ahler manifold $Q$, and the tangent space at a solution $(T_0,T_1,T_2,T_3)$ can be identified with the space of solutions to the following system of linear
equations:
\begin{equation}\begin{array}{c} \dot{t}_0+[T_0,t_0]+[T_1,t_1]+[T_2,t_2]+[T_3,t_3]=0,\\
\dot{t}_1+[T_0,t_1]-[T_1,t_0]-[T_2,t_3]+[T_3,t_2]=0,\\
\dot{t}_2+[T_0,t_2]+[T_1,t_3]-[T_2,t_0]-[T_3,t_1]=0,\\
\dot{t}_3+[T_0,t_3]-[T_1,t_2]+[T_2,t_1]-[T_3,t_0]=0.\end{array}\label{tangent}\end{equation}
The first equation is the condition that $(t_0,t_1,t_2,t_3)$ is orthogonal to the infinitesimal gauge transformations and the remaining three are linearisations of Nahm's equations.

\subsection{Curvature of the Nahm moduli spaces}
We shall estimate the sectional curvature of the above moduli space $Q$ of solutions to Nahm's equations from Theorem \ref{main}. We need first some definitions.
\par
\begin{definition} Let $\lambda:[a,b]\rightarrow \oR$ be a continuous function. Let $E=\{u\in W^{2,1}(a,b); u(a)=u(b)=0\}$ and let
$L_\lambda:E\rightarrow L^1(a,b)$ be the linear operator
$$L_\lambda(u)(s)=\ddot{u}(s) -\lambda(s)^2u(s).$$
$L$ has a continuous inverse and we define $N(\lambda)$ as the norm of $j\circ L_\lambda^{-1}$ where $j:W^{2,1}(a,b)\rightarrow C^0([a,b])$ is
the embedding (and $C^0([a,b])$ is equipped with usual $\max$ norm).\label{N}\end{definition}

\begin{theorem} Let $T=(T_0,T_1,T_2,T_3)$ be a solution to Nahm's equations on $[a,b]$ which is an element of a hyperk\"ahler
quotient $Q$ constructed above. Let $\lambda(s)$ be a continuous real function such that, for every $s\in [a,b]$, $\lambda^2(s)$ is not
greater than the smallest eigenvalue of the operator  $H(s)=-(\ad T_1(s))^2 -(\ad T_2(s))^2-(\ad T_3(s))^2$. Then the sectional curvature
of $Q$ at $T$ is bounded by $18N(\lambda)^{1/2}$.\label{T}
\end{theorem}
\begin{proof}
We shall use Propositions \ref{main},\ref{main2} together with Lemmae \ref{H} and \ref{estimate}. For every solution $T$ to Nahm's equations the norm $\|\cdot\|_T$ on $\Lie \sG$ will be
the $L^\infty$-norm, so the dual norm $\|\cdot\|_T^\ast$ is the $L^1$-norm.
Since the Levi-Civita connection on $\sA$ is simply the directional derivative, we have, using \eqref{check}, for a horizontal tangent vector
$X=(t_0,t_1,t_2,t_3)$,
$$\nabla_X\check{\rho}=\bigl([\rho,t_0],[\rho,t_1],[\rho,t_2],[\rho,t_3]\bigr).$$
Hence, because of the normalisation \eqref{normalise},
\begin{equation} g\bigl(\nabla_X\check{\rho},\nabla_X\check{\rho}\bigr)\leq 4g(X,X)\|\rho\|_{L^\infty}^2,\label{nabla}\end{equation}
and consequently
\begin{equation} F(T)\leq 2.\label{2}\end{equation}
We now estimate $l(T)$ using the second statement in Lemma \ref{estimate}. We first observe that under the action of a gauge transformation
$g(s)$, $\rho,\check{\rho}$ and $\Lambda(\check{\rho})$ are all pointwise conjugated by $g(s)$ and therefore the values of $l$ at $T$ and at
$g.T$ are the same.  Similarly,  the eigenvalues of $H$ do not change under the action of $g$. Therefore we can assume that $T_0\equiv 0$.
\par
The form $\Lambda$ is given by the  equation \eqref{perp}.
Substituting \eqref{check} into \eqref{perp}   we get (as $T_0=0$)
$$ \ddot{\rho}=H(s)\rho -\Lambda(\check{\rho}),$$
where $H(s)$ is the positive-definite Hermitian operator defined in the statement of the theorem. Therefore
$$\frac{d^2}{ds^2}\langle \rho,\rho\rangle=2\langle \dot{\rho},\dot{\rho}\rangle +2 \langle \ddot{\rho},\rho\rangle\geq
2\langle \dot{\rho}(s),\dot{\rho}(s)\rangle+2\lambda(s)^2|\rho(s)|^2-|\Lambda(\check{\rho})(s)||\rho(s)|.$$
Combining this with
$$\frac{d^2}{ds^2}|\rho(s)|^2=2|\rho(s)|\frac{d^2|\rho(s)|}{ds^2}+2\left(\frac{d|\rho(s)|}{ds}\right)^2\leq
2|\rho(s)|\frac{d^2 |\rho(s)|}{ds^2} +2\left|\frac{d\rho(s)}{ds}\right|^2,$$
we get
$$\frac{d^2}{ds^2}|\rho(s)|\geq \lambda(s)^2|\rho(s)|-|\Lambda(\check{\rho})(s)|.$$

From this one easily concludes that $|\rho(s)|$ is point-wise bounded by the solution $u(s)$ to $L_\lambda(u)=-|\Lambda(\check{\rho})(s)|$ with $u(a)=u(b)=0$. The definition of $N(\lambda)$ gives now $\|\rho\|_{L^\infty}\leq
N(\lambda)\|\Lambda(\check{\rho})\|_{L^1}$, and Lemma \ref{estimate} implies that $l(T)\leq N(\lambda)^{1/2}$.
\end{proof}

\begin{corollary} The sectional curvature of a moduli space of solutions to Nahm's equations on $(a,b)$ with the boundary conditions
 described above and the metric \eqref{metric} satisfying \eqref{normalise} is bounded by $9\sqrt{b-a}$.
\end{corollary}
\begin{proof} We use the last theorem with $\lambda=0$. Thus we need to estimate the $C^0$-norm of the unique solution $u$ of
$\ddot{u}(s)=h(s)$ with $u(a)=u(b)=0$ in terms of the $L^1$-norm of $h$. We have an explicit solution to this boundary problem:
$$ u(s)=-\frac{1}{b-a}(b-s)\int_a^s(\tau-a)h(\tau)d\tau-\frac{1}{b-a}(s-a)\int_s^b(b-\tau)h(\tau)d\tau.$$
From this
$$|u(s)|\leq \frac{(b-s)(s-a)}{b-a}\|h\|_{L^1}$$
and so $N(0)=(b-a)/4$.
\end{proof}

\begin{remarknon} The moduli space of charge $k$ $SU(2)$-monopoles arises when $G=U(k)$ and the residues $\alpha_i,\beta_i$ of the solutions to Nahm's
equations at $a,b$ define irreducible representations of $\su(2)$ \cite{Nahm}. Therefore the curvature of this moduli space is bounded. One
can extend the above proof to more complicated moduli spaces of solutions to Nahm's equations, such as those in \cite{Hur}, to show that
the moduli spaces of $G$-monopoles with maximal symmetry breaking have bounded curvature ($G$ - a classical compact group).\end{remarknon}

\begin{remarknon} If we are only interested in a global bound on the curvature, and not in the finer estimates of Theorem \ref{T}, then the proof is much simpler. We can use Corollary \ref{O'Neill}: if $|\check{\rho}_T|\leq 1$, then $\|\dot{\rho}\|_{L^2}\leq 1$ and now the Poincar\'e inequality together with the Sobolev embedding $W^{1,2}(a,b)\rightarrow C^0([a,b])$ implies that $\|\rho\|_{L^\infty}\leq K$ for some constant $K$ depending only on $b-a$. The estimate \eqref{nabla} and Corollary \ref{O'Neill} give us a bound on the sectional curvature.\end{remarknon}

\section{Moduli spaces of solutions to Hitchin's equations}

We briefly recall the hyperk\"ahler quotient construction of the moduli space of solutions to Hitchin's equations on a Riemann surface with  boundary \cite{Hit,Don2}. Since the fields do not have singularities, this is actually simpler than for Nahm's equations considered in the previous section.

Let $Z$ be a compact connected $2$-dimensional K\"ahler manifold $Z$ with  boundary $\partial Z$ which may be empty. Let $P$ be a principal bundle over $Z$ with a compact structure group $G$ and a $G$-trivialisation over $\partial Z$. For a positive integer $k$, we denote by $\Omega^1_k(\ad P)$ the space of $\ad P$-valued $1$-forms of Sobolev class $W^{k,2}$, i.e. those whose first $k$ derivatives are square-integrable. Similarly for other forms. We now consider the affine manifold $\sM=\bigl(D_0+\Omega^{0,1}_4(\ad P\otimes \oC)\bigr)\oplus \Omega^{0,1}_4(\ad P\otimes \oC)$, where $D_0$ is a fixed $G$-connection. We identify $\Omega_4^{0,1}\bigl(\ad P\otimes \oC)$ with $\Omega_4^1(\ad P)$ and we view elements of  $D_0+\Omega^{0,1}_4(\ad P\otimes \oC)$ as $G$-connections. The metric on $Z$ and an invariant inner product on $\fG$ induce an $L^2$-metric on $\Omega_4^{1,0}(\Ad P\otimes \oC)$ and on $\Omega_4^{0,1}(\Ad P\otimes \oC)$. With
this metric, $\sM$ is a flat weak hyperk\"ahler Hilbert manifold. The three anti-commuting complex structures  are given, on each tangent space $\Omega^{0,1}\bigl(\Ad P\otimes \oC\bigr)\oplus \Omega^{1,0}\bigl(\Ad P\otimes \oC\bigr)$, by:
$$ I_1(a,\phi)=(ia,i\phi),\enskip I_2(a,\phi)=(-\phi^\ast,a^\ast), \enskip I_3(a,\phi)=(-i\phi^\ast,ia^\ast).$$
The Levi-Civita connection is again provided by the directional derivative.
\par
The gauge group $\sG$ consists of $G$-valued gauge transformations of Sobolev class $W^{5,2}$ which are identity on $\partial Z$. The action of this group is  isometric and tri-Hamiltonian with  hyperk\"ahler moment map $\mu=(\mu_1,\mu_2,\mu_3)$ described by \cite{Hit, Don2}, as giveb by the Hitchin equations, i.e.
\begin{equation}  \mu_1\bigl(A,\Phi\bigr)= F_A+\bigl[\Phi,\Phi^\ast\bigr],\quad(\mu_2+i\mu_3)\bigl(A,\Phi\bigr)=\bar{\partial}_A\Phi,\label{Hitchin}\end{equation}
where $F_A$ is the curvature of $A$. The action is proper and clearly free if $Z$ is connected with non-empty boundary (since the elements of $\sG$ are identity on $\partial Z$). If $\partial Z=\emptyset$, then the action is free apart from reducible pairs $(A,\Phi)$. The same restriction holds for the ellipticity of the action. Indeed, the operator $\rho\mapsto \Lambda(\check{\rho}_m)$ considered in \eqref{Ph} is identified, at $m=(A,\Phi)$ with
$$ \rho\mapsto i\bar{\partial}_A\partial_A\rho+\bigl[\Phi,[\Phi^\ast, \rho]\bigr].$$
 This operator is self-adjoint and positive apart from the case when $\partial Z=\emptyset$ and $(A,\Phi)$ is reducible. Thus, away from such points $m$, we get a unique solution to $\Lambda(\check{\rho}_m)=h\in \Omega^{1,1}_3(\Ad P\otimes \oC)$ with $\rho\in \Lie \sG$ and so the action of $\sG$ is elliptic. Propositions  \ref{stronger} and \ref{quot} show now that that the quotient $Q=\mu^{-1}(0)/\sG$, where $\mu$ is given by \eqref{Hitchin}, is a weak hyperk\"ahler manifold away from the set points where $\sG$ has a non-trivial isotropy. In the case of a closed Riemann surface, $Q$ is finite-dimensional \cite{Hit}, while for the moduli space on (connected) $Z$ framed on $\partial Z$ it is a smooth infinite-dimensional manifold  \cite{Don2}. In the latter case, when $Z$ is the closed disc, Donaldson  shows that $Q$ can be identified with a complex bundle over the based loop space $\Omega G$. For a general $Z$, one would expect that $Q$ is a complex vector bundle over the space $\text{Maps}(\partial Z,G^\cx)/\text{Hol}(Z,G^\cx)$, the latter being a submanifold of $Q$ obtained by setting $\Phi=0$ \cite{Don2}.
\par
We now have the following analogue of Theorem \ref{T}:
\begin{theorem} Let $(A,\Phi)$ be a solution to Hitchin's equations on $Z$. The sectional curvature of the moduli space $\mu^{-1}(0)/\sG$ at a point $m=(A,\Phi)$ is bounded by
$$ C\cdot\sup_{\rho\in \Lie\sG- \{0\}}\frac{\|\rho\|_{W^{1,2}}}{|\check{\rho}_m|}= C\cdot\sup\Bigl\{\|\rho\|_{W^{1,2}}\enskip;\enskip   \bigl\|[\rho,A]\bigr\|^2_{L^{2}}+ \bigl\|[\rho,\Phi]\bigr\|^2_{L^{2}}=1\Bigr\},$$
where $C$ does not depend on $(A,\Phi)$.
\end{theorem}
{\em Proof.} According to Corollary \ref{O'Neill}, we need to prove that the $L^2$-norm of $\nabla_X\check{\rho}$ is bounded at $(A,\Phi)$ by $C\|\rho\|_{W^{1,2}}$  for any $X,\check{\rho}$ whose $L^2$-norms are $1$. Since $\nabla$ is simply the directional derivative, it follows that for $X=(a,\phi)$
$$\nabla_X\check{\rho}=\bigl([\rho,a],[\rho,\phi]\bigr).$$ The desired bound on $\nabla_X\check{\rho}$ follows, by using local trivialisations and partitions of unity, from the following
extension of the Sobolev multiplication theorem to the critical case:

\begin{proposition} Let  $Y$ be a compact $2$-dimensional Riemannian manifold $Y$ with or without boundary. Then the multiplication of functions is a well-defined  continuous bilinear operator
$$ L^2(Y)\times W^{1,2}(Y)\rightarrow L^2(Y).$$
\end{proposition}
\begin{proof} Although we do not have an embedding of $W^{1,2}(Y)$ into $L^\infty(Y)$, we do have a  continuous embedding into certain Orlicz space, namely, if $\psi\in W^{1,2}(Y)$, then $e^{|\psi|}\in L^1(Y)$ and
\begin{equation} \int_Y e^{|\psi|}d\mu\leq C\exp\bigl(\alpha\|\psi\|_{W^{1,2}}^2\bigr)\label{exp} \end{equation}
for some constants $C,\alpha$ depending only on $Y$
(see, e.g., Theorem 2.46 in \cite{Aubin} for a proof). Let now $\phi\in L^2(Y)$ and  $\psi\in W^{1,2}(Y)$ both have norm $1$ in the respective spaces. Define, for a nonnegative integer $K$,
$$ Y_K=\bigl\{z\in Y;\enskip K\leq |\psi(z)|\leq K+1\bigr\}.$$
It follows from \eqref{exp} that $\mu(Y_K)\leq Ce^{\alpha -K}$. We now compute
$$\int_Z|\phi\psi|^2d\mu=\sum_{K\geq 0} \int_{Y_K}|\phi\psi|^2d\mu \leq \sum_{K\geq 0} (K+1)^2Ce^{\alpha -K}$$
and, hence, the $L^2$-norm of $\phi\psi$ is bounded by some constant depending only on $Y$.
\end{proof}

\begin{corollary} Let $Z$ be a compact connected Riemann surface with a non-empty boundary and $P$ a principal $G$-bundle trivialised on $\partial Z$.
Then the sectional curvature of the moduli space $Q$ of solutions to Hitchin's equations on $(Z,\Ad P\otimes \oC)$ framed on $\partial Z$ is bounded.\end{corollary}
\begin{proof} Let $(A,\Phi)\in Q$ and $\check{\rho}$ be a fundamental vector field whose $L^2$-norm is $1$ at  $(A,\Phi)$. Then $\|\nabla_A\rho\|_{L^2}\leq 1$ and, since  $\rho$ vanishes on the boundary, the Kato  and  Poincar\'e inequalities imply that the $W^{1,2}$-norm of $\rho$ is bounded by some constant, which does not depend on $(A,\Phi)$. Therefore the expression in the statement of the last theorem is finite and does not depend on $(A,\Phi)$. \end{proof}

\begin{ack} The author thanks the Humboldt Foundation for a Fellowship, during which this paper was partly written.\end{ack}


\begin{thebibliography}{99}


\bibitem{Aubin}
{T. Aubin}, {\em Nonlinear analysis on manifolds. Monge-Amp\`ere equation},  Springer Verlag, Berlin (1982).

\bibitem{Besse}
{A.L. Besse}, {\em Einstein manifolds}, Springer Verlag, Berlin (1987).


\bibitem{Don}{ S.K. Donaldson}, `Nahm's equations and the classification of monopoles', {\it Comm. Math. Phys.} 96 (1984), 387--407.

\bibitem{Don2}{ S.K. Donaldson}, `Boundary value problems for Yang-Mills fields', {\it J. Geom. Phys.}  8  (1992),  89--122.

\bibitem{Hit}{N.J. Hitchin}, ` The self-duality equations on a Riemann surface',
{\it Proc. London Math. Soc.} (3) 55 (1987),  59--126.

\bibitem{Hur}
{ J.C. Hurtubise}, `The classification of monopoles for the classical groups', {\it Comm. Math. Phys.} 120 (1989), 613--641.

\bibitem{JP}
{J. Jost and X.-W. Peng}, `Group actions, gauge transformations, and the calculus of variations', {\it Math. Ann.} 293 (1992),  595--621.

\bibitem{L}
{S. Lang}, {\em Fundamentals of differential geometry}, Springer Verlag, Berlin (1999).

\bibitem{Nahm}
{ W. Nahm}, `The construction of all self-dual monopoles by the ADHM method', in {\em Monopoles in quantum field theory}, World Scientific,
Singapore (1982).

\bibitem{Swa}
{C. Swartz}, `Continuity and hypocontinuity for bilinear maps', {\it. Math. Z.} 186 (1984), 321--329.


\end{thebibliography}
\end{document}